\documentclass[preprint]{elsarticle}
\usepackage{amsbsy}
\usepackage{amsfonts}
\usepackage{amsmath}
\usepackage{amssymb}

\newdefinition{dfn}{Definition}
\newtheorem{theorem}{Theorem}
\newtheorem{stw}{Proposition}

\newproof{dwd}{Proof}
\newtheorem{wn}{Corollary}
\newtheorem{lem}{Lemma}
\newdefinition{uw}{Remark}
\newdefinition{prob}{Problem}

\newdefinition{prk}{Example}

\begin{document}

\begin{frontmatter}
\title{On the proportion of prefix codes in the set of   three-element codes}
\author[rvt]{Adam Woryna}
\ead{adam.woryna@polsl.pl}
\address[rvt]{Silesian University of Technology, Faculty of Applied Mathematic, ul. Kaszubska 23, 44-100 Gliwice, Poland}

\begin{abstract}
Let $L$ be a finite sequence of natural numbers. In \cite{4,22}, we derived some interesting properties for the ratio $\rho_{n,L}=|PR_n(L)|/|UD_n(L)|$, where $UD_n(L)$ denotes the set of all codes over an $n$-letter alphabet and  with length distribution $L$, and  $PR_n(L)\subseteq UD_n(L)$ is the corresponding subset of  prefix codes. In the present paper, we study the case when the length distributions are three-element sequences. We show  in this case that the ratio $\rho_{n,L}$ is always greater than $\alpha_n$, where $\alpha_n=(n-2)/n$ for $n>2$ and $\alpha_2=1/6$. Moreover, the number $\alpha_n$ is the best possible lower bound for this ratio, as  the length distributions of the form $L=(1,1,c)$ and $L=(1,2,c)$ assure that  the  ratios asymptotically approach $\alpha_n$. Namely,  if $L=(1,1,c)$, then  $\rho_{n,L}$ tends to  $(n-2)/n$   with $c\to\infty$, and,  if $L=(1,2,c)$, then $\rho_{2,L}$ tends to $1/6$ with $c\to\infty$.
\end{abstract}

\end{frontmatter}

\section{Motivation and the results}

A  {\it code}  over a finite alphabet $X$ is a finite sequence $C=(v_1,\ldots, v_m)$  of words over $X$ (so-called {\it  code-words}) such that every  $w\in X^*$ has at most one factorization into the code-words, i.e. if $w=v_{i_1}\ldots v_{i_l}=v_{j_1}\ldots v_{j_{l'}}$ for some $l,l'\geq 1$ and $1\leq i_t,j_{t'}\leq m$ ($1\leq t\leq l$, $1\leq t'\leq l'$), then $l'=l$ and $i_t=j_t$ for every $1\leq t\leq l$ (note that our definition differs a bit from the more usual one, where codes are considered as the sets of words rather than the sequences -- see also~\cite{17,0}). A code $C=(v_1,\ldots, v_m)$ is called a {\it prefix code} if it satisfies the following condition: for all $1\leq i,j\leq m$  the code-word $v_i$ is a prefix (initial segment) of the code-word $v_j$  if and only if $i=j$.

For every natural number $n\geq 2$ and every finite sequence $L=(a_1, \ldots, a_m)$ of natural numbers, we consider the set $UD_n(L)$ of all codes  over an $n$-letter alphabet $X$ with length distribution  $L$, i.e. a code  $C=(v_1, \ldots, v_m)$ over $X$ belongs to $UD_n(L)$  if and only if $|v_i|=a_i$ for every $1\leq i\leq m$. We denote by $PR_n(L)$ the corresponding subset of prefix codes in the set $UD_n(L)$.

The prefix codes form the most useful and important class of codes. Therefore it is natural to ask about the contribution of these types of  codes in various classes of codes. For  given $L$ and $n\geq 2$ both the sets $UD_n(L)$ and $PR_n(L)$ are finite, and  such a contribution may be defined  as the ratio
$$
\rho_{n,L}=\frac{|PR_n(L)|}{|UD_n(L)|}.
$$
Note that if $L=(a_1, \ldots, a_m)$, then by McMillan theorem (\cite{6}), the following statements are equivalent:
\begin{itemize}
\item $UD_n(L)\neq \emptyset$.
\item $PR_n(L)\neq \emptyset$.
\item $\sum_{i=1}^mn^{-a_i}\leq 1$ (the so-called Kraft inequality).
\end{itemize}
If $L$ is constant, then the code-words of every  code $C=(v_1, \ldots, v_m)$ in $UD_n(L)$  have  the same length. In particular,
for all $1\leq i,j\leq m$, the code-word $v_i$ is a prefix of the code-word $v_j$ if and only if $v_i=v_j$, which implies that  $C$ is also a prefix code. Thus, if $L$ is constant, the equality $UD_n(L)=PR_n(L)$ holds. It follows that for every $n\geq 2$ and $m\geq 1$, there is a sequence $L$ of length $|L|=m$ such that  the ratio $\rho_{n,L}$ is equal to $1$. On the other hand, as we showed in~\cite{4}, the following theorem holds:
\begin{theorem}[\cite{4}, Theorem~1]\label{twt1}
If  $L$ is non-constant and $UD_n(L)\neq\emptyset$, then
$$
\frac{1}{\rho_{n,L}}\geq 1+\frac{r_a\cdot r_b}{n^{a+b}-n^{\max\{a,b\}}},
$$
where $a$ and $b$ are arbitrary two different values of $L$ and $r_a$ (resp. $r_b$) is the number of those elements in $L$ which are equal to $a$ (resp. to $b$).
\end{theorem}

Obviously, for every $n\geq 2$ and $m\geq 1$ there are infinitely many sequences $L$ of length $m$ such that the sets $UD_n(L)$ and $PR_n(L)$ are non-empty. Hence, given $n$ and $m$, it is  non-trivial to ask whether $\rho_{n,L}$ can be arbitrarily close to zero (depending on $L$). In the paper~\cite{22}, we negatively answered this question by  the following result:
\begin{theorem}[\cite{22}, Theorem 2]
If $n\geq 2$ and  $m\geq 1$, then for every  sequence $L$ of length $m$ such that $UD_n(L)\neq \emptyset$, the following inequality holds
$$
\rho_{n,L}\geq q_{n,m}\cdot \left(\frac{n-(m)_{n-1}}{n^{\left\lfloor \frac{m}{n-1}\right\rfloor+1}}\right)^{m-1},
$$
where $(m)_{n-1}$ is the remainder from the division of $m$ by $n-1$ and
$$
q_{n,m}:=\left\{
\begin{array}{ll}
1,&n\geq m,\\
\frac{(m-1)!}{(m-1)^{m-1}},&n<m.
\end{array}
\right.
$$
\end{theorem}
By the above theorem, we see that for every $n\geq 2$ and $m\geq 1$ the infimum
\begin{equation}\label{infim}
\inf_{L}\rho_{n,L}
\end{equation}
taken over all  sequences $L$ of length $m$ such that $UD_n(L)\neq \emptyset$ is greater than zero (note that this infimum  depends only on $n$ and $m$). In the paper~\cite{22}, we derived various interesting properties for this infimum. For example, we showed there that for every $n\geq 2$, it tends to $0$ when $m\to\infty$, and for every $m\geq 1$, it tends to $1$ when $n\to\infty$.

In~\cite{22}, we also derived for all $a,b\geq 1$ the equality
\begin{equation}\label{f1}
|UD_n((a,b))|=n^{a+b}-n^{{\rm gcd}(a,b)},
\end{equation}
which we apply in the proof of the following
\begin{theorem}[\cite{22}, Corollary~5]
For every $n\geq 2$ the infimum (\ref{infim}) taken over all sequences $L$ of length $|L|=2$ for which $UD_n(L)\neq \emptyset$ is equal to $(n-1)/n$.
\end{theorem}

The formula  (\ref{f1})   follows from a nice characterization of two-element codes. Namely, a sequence $C=(w,v)$ is a code if and  only if $vw\neq wv$, or, equivalently, the words $w$, $v$ are not the powers of the same word (see also~\cite{17} or~\cite{1}).
The situation is much more complicated in the case  of codes of length three, as the full  characterization of such  codes is not known (for some partial results see~\cite{14,9,11}). Nevertheless, the properties of  three-element codes are studied within the years (\cite{12,14,2,8,10,16}) and  some problems remain still open for these codes (see~\cite{15,16,13}).

For every $n\geq 2$ let us define the number $\alpha_n$ as follows
$$
\alpha_n:=\left\{
\begin{array}{ll}
1/6,&{\rm if}\;n=2,\\
(n-2)/n, &{\rm if}\;n>2.
\end{array}
\right.
$$
For the main results of the present paper, we prove the following theorems
\begin{theorem}\label{t1}
If  $n\geq 2$ and  $L$ is an arbitrary sequence of length   three such that $UD_n(L)\neq\emptyset$, then $\rho_{n,L}>\alpha_n$.
\end{theorem}

\begin{theorem}\label{t2}
For every $c\geq 1$ let us denote $L_c=(1,1,c)$. Then for every $n>2$, we have:
$$
\lim_{c\to\infty}\rho_{n,L_c}=\frac{n-2}{n}=\alpha_n.
$$
\end{theorem}

\begin{theorem}\label{t3}
For every $c\geq 2$ let us denote $L_c=(1,2,c)$. Then  we have:
$$
\lim_{c\to\infty}\rho_{2,L_c}=\frac{1}{6}=\alpha_2.
$$
\end{theorem}
As a direct consequence of Theorems~\ref{t1}-\ref{t3}, we obtain the following
\begin{wn}
For every $n\geq 2$ the infimum (\ref{infim}) taken over all sequences $L$ of length $|L|=3$ for which $UD_n(L)\neq \emptyset$ is equal to $\alpha_n$.
\end{wn}

\section{The proof of Theorems~\ref{t1}-\ref{t2}}\label{sec2}

In this section, we give the proof of Theorems~\ref{t1}--\ref{t2}. To this aim,  we need at first the following two propositions.

\begin{stw}\label{prefi}
If $n\geq 2$ and $L=(a,b,c)$ is an arbitrary three-element sequence such that $PR_n(L)\neq\emptyset$, then
$$
|PR_n(L)|=n^c(n^{a+b}-2n^{b}-n^{a}+n^{b-a}+1).
$$
\end{stw}
\begin{dwd}[of Proposition~\ref{prefi}]
Since $|PR_n(L')|=|PR_n(L)|$ for every  sequence  $L'$ obtained from $L$ by  permuting the elements, we can  assume that $a\leq b\leq c$. An arbitrary prefix code $(w,v,u)\in PR_n(L)$ can be constructed as follows. The word $w$ can be freely chosen among the $n^a$ words of length $a$. The word $v$ can be freely chosen among the words of length $b$ which do not have $w$ as a prefix. The number of these words is equal to $n^b-n^{b-a}$. Finally, the word $u$ can be freely chosen among the words of length $c$ which do not have any of the words $w,v$ as a prefix. Since the number of words of length $c$ which have one of the words $w, v$ as a prefix is equal to $n^{c-a}+n^{c-b}$, we can choose  $u$  among the $n^c-n^{c-a}-n^{c-b}$ available words. In consequence, we obtain
$$
|PR_n(L)|=n^a\left(n^b-n^{b-a}\right)\left(n^c-n^{c-a}-n^{c-b}\right),
$$
and the desired formula follows  by  multiplying the above brackets.\qed
\end{dwd}

\begin{stw}\label{l11}
$|UD_n((1,1,c))|=n(n-1)(n^c-2^c)$ for every $n\geq 2$ and $c\geq 1$.
\end{stw}
\begin{dwd}[of Proposition~\ref{l11}]
The claim  follows from the observation that a sequence $(x,y,w)\in X\times X\times X^c$ is a code if and only if $x\neq y$ and $w\notin \{x,y\}^c$. Obviously, if at least one of these  conditions does not hold, then this sequence is not a code. Conversely, let us assume that $x\neq y$ and $w\notin \{x,y\}^c$. To show that $(x,y,w)\in UD_n((1,1,c))$, we use  the Sardinas-Patterson algorithm (\cite{3}). Namely,  we define $\mathcal{D}_0:=\{x,y,w\}$ and for each $i\geq 1$ we define $\mathcal{D}_i$ as the set of all non-empty words $u\in X^*$ for which  the following condition holds: $\mathcal{D}_{i-1}u\cap \mathcal{D}_0\neq \emptyset$ or $\mathcal{D}_0u\cap \mathcal{D}_{i-1}\neq\emptyset$, where $\mathcal{D}_iu:=\{vu\colon v\in \mathcal{D}_i\}$ for $i\geq 0$. According to the Sardinas-Patterson theorem, we have $(x,y,w)\in UD_n((1,1,c))$ if and only if $\mathcal{D}_i\cap \mathcal{D}_0=\emptyset$ for every $i\geq 1$. Since $w\notin \{x,y\}^c$, there is the smallest number   $1\leq i_0\leq c$  such that the $i_0$-th letter of $w$ does not belong to the set $\{x,y\}$. For each $1\leq i<i_0$ let $w_i$ be the suffix of $w$ of length $c-i$. Because of the minimality of  $i_0$, none of the words $w_i$ ($1\leq i<i_0$)  is a prefix of $w$. Hence,  by the trivial induction on $i$, we obtain  $\mathcal{D}_i=\{w_i\}$ for every $1\leq i<i_0$ and $\mathcal{D}_i=\emptyset$ for every $i\geq i_0$. Thus  $\mathcal{D}_i\cap \mathcal{D}_0=\emptyset$ for each $i\geq 1$, and hence $(x,y,w)\in UD_n((1,1,c))$. \qed
\end{dwd}

We are ready now to prove our main results.
{
\renewcommand{\thetheorem}{\ref{t1}}
\begin{theorem}
If  $n\geq 2$ and  $L$ is an arbitrary sequence of length   three such that $UD_n(L)\neq\emptyset$, then $\rho_{n,L}>\alpha_n$.
\end{theorem}
\addtocounter{theorem}{-1}
}
\begin{dwd}
Let $L=(a,b,c)$. Since $\rho_{n, L}=\rho_{n,L'}$ for  every  sequence  $L'$ obtained from $L$ by  permuting the elements, we can  assume that $a\leq b\leq c$. For every $(w,v,u)\in UD_n(L)$, we have $(w,v)\in UD_n((a,b))$, and hence, by the formulae~(\ref{f1}), we obtain
\begin{equation}\label{w4}
|UD_n(L)|\leq n^c\cdot |UD_n((a,b))|=n^c(n^{a+b}-n^{{\rm gcd}(a,b)})\leq n^c(n^{a+b}-n).
\end{equation}
By Proposition~\ref{prefi} and by the inequality~(\ref{w4}), we obtain:
$$
\rho_{n,L}=\frac{|PR_n(L)|}{|UD_n(L)|}\geq \frac{n^{a+b}-2n^{b}-n^{a}+n^{b-a}+1}{n^{a+b}-n}.
$$

Let us denote by $Q(a,b)$ the quotient on the right side of the above inequality. Then we have
\begin{eqnarray*}
Q(a,b)&=&\frac{n^{a+b}\left(1-2n^{-a}+n^{-2a}\right)-n^{a}+1}{n^{a+b}-n}=\\
&=&\frac{\left(n^{a+b}-n+n\right)\left(1-2n^{-a}+n^{-2a}\right)-n^{a}+1}{n^{a+b}-n}=\\
&=&\frac{\left(n^{a+b}-n\right)\left(1-2n^{-a}+n^{-2a}\right)+n\left(1-2n^{-a}+n^{-2a}\right)-n^{a}+1}{n^{a+b}-n}=\\
&=&\left(1-2n^{-a}+n^{-2a}\right)+\frac{R(a)}{n^{a+b}-n}=\left(1-\frac{1}{n^a}\right)^2+\frac{R(a)}{n^{a+b}-n},
\end{eqnarray*}
where
\begin{eqnarray*}
R(a):=n\left(1-2n^{-a}+n^{-2a}\right)-n^a+1=n\left(1-\frac{1}{n^a}\right)^2-n^a+1=\\
=\frac{n\left(n^a-1\right)^2}{n^{2a}}-\left(n^a-1\right)=\frac{n^a-1}{n^{2a-1}}\left(n^a-1-n^{2a-1}\right)<0.
\end{eqnarray*}
Since $R(a)<0$, it follows by the equality
$$
Q(a,b)=\left(1-\frac{1}{n^a}\right)^2+\frac{R(a)}{n^{a+b}-n}
$$
that $Q(a,b)$ increases with $b$.  Hence, if $b\geq a+1>2$, then
\begin{eqnarray*}
\rho_{n,L}&\geq& Q(a,b)\geq Q(a,a+1)=\frac{n^{2a+1}-2n^{a+1}+n^{a+1-a}-n^{a}+1}{n^{2a+1}-n}=\\
&=&1-\frac{2n+1}{n^{a+1}+n}>1-\frac{2n+1}{n^2+n}\geq \max\left\{\frac{n-2}{n},\frac{1}{6}\right\}\geq \alpha_n.
\end{eqnarray*}
If $b\geq a+1=2$, then for every $(w,v,u)\in UD_n(L)$ we have $u\neq w^c\in X^c$, and hence
$$
|UD_n(L)|\leq (n^c-1)\cdot |UD_n((a,b))|<n^c(n^{a+b}-n).
$$
Consequently, we have in this case
\begin{eqnarray*}
\rho_{n,L}&=&\frac{|PR_n(L)|}{|UD_n(L)|}>Q(a,b)\geq Q(a,a+1)=Q(1,2)=\\
&=&1-\frac{2n+1}{n^2+n}\geq \max\left\{\frac{n-2}{n},\frac{1}{6}\right\}\geq \alpha_n.
\end{eqnarray*}
If $b=a>1$, then
\begin{eqnarray*}
\rho_{n,L}\geq Q(a,b)=Q(a,a)=\frac{(n^a-1)(n^a-2)}{n^{2a}-n}>\frac{(n^a-1)(n^a-2)}{n^{2a}-1}=\\
=\frac{n^a-2}{n^a+1}=1-\frac{3}{n^a+1}\geq 1-\frac{3}{n^2+1}>\max\left\{\frac{n-2}{n},\frac{1}{6}\right\}\geq \alpha_n.
\end{eqnarray*}
If $b=a=1$, then $n>2$ and by Proposition~\ref{l11}, we have
\begin{equation}\label{ff}
\rho_{n,L}=\rho_{n,(1,1,c)}=\frac{n^c(n^{2}-3n+2)}{n(n-1)(n^c-2^c)}=\frac{n-2}{n}\cdot\frac{1}{1-(2/n)^c}>\alpha_n.
\end{equation}
Thus in every case, we have $\rho_{n,L}>\alpha_n$, which finishes the proof of Theorem~\ref{t1}.\qed
\end{dwd}

{
\renewcommand{\thetheorem}{\ref{t2}}
\begin{theorem}
For every $c\geq 1$ let us denote $L_c=(1,1,c)$. Then for every $n>2$, we have:
$$
\lim_{c\to\infty}\rho_{n,L_c}=\frac{n-2}{n}=\alpha_n.
$$
\end{theorem}
\addtocounter{theorem}{-1}
}
\begin{dwd}
The proof is immediate, as we have by~(\ref{ff}):
$$
\lim\limits_{c\to\infty}\rho_{n,L_c}=\lim\limits_{c\to\infty}\frac{n-2}{n(1-(2/n)^c)}=\alpha_n.
$$\qed
\end{dwd}

\section{The proof of Theorem~\ref{t3}}\label{sec3}

In this section, we derive Theorem~\ref{t3} describing the case of  binary codes. Surprisingly, this is the most involving step of our study. As we will see, the main burden of the proof relies on some two propositions below (Proposition~\ref{lll3} and Proposition~\ref{lll6}) deriving  the cardinalities of two specially constructed subsets. In the proof of both these propositions, we will use the following auxiliary  lemma.

\begin{lem}\label{lll}
Let $C=(v_1, \ldots, v_m)$ be a sequence of non-empty words with the following property: There exist $1\leq \mu,\kappa\leq m$, $\mu\neq \kappa$ such that $v_\mu=v_\kappa v$ for some $v\in X^*$. If $C$ is not a code, then the sequence $C':=(v_1', \ldots, v_m')$ is not a code, where $v'_i:=v_i$ for $i\neq \mu$ and $v'_\mu:=v$.
\end{lem}
\begin{dwd}[of Lemma~\ref{lll}]
Obviously, if any word $w$ has a factorization into  words from the sequence  $C$, then it has also a factorization into  words from the sequence $C'$. In fact, if $w=w_1\ldots w_l$, where all $w_i$ ($1\leq i\leq l$) are from $C$ and if $1\leq i_0\leq l$ is the smallest number such that $w_{i_0}=v_\mu=v_\kappa v$, then there exists $l'>l$ and the  words $W_1, W_2, \ldots, W_{l'}$ in $C'$  such that $W_i=w_i$ for every $1\leq i<i_0$, $W_{i_0}=v_\kappa$,  $W_{i_0+1}=v$ and $w_1\ldots w_l=W_1\ldots W_{l'}$.

Suppose contrary that $C$ is not a code and $C'$ is a code. Since   each word in $C$ is non-empty, we have $v\neq v_\mu$. This implies that the word  $v$ does not belong to $C$ (otherwise  $v=v_i=v'_i$ for some $i\neq \mu$,  and then $C'$ would not be a code, contrary to our assumption).

Let $w$ be the shortest word with  the following two different factorizations
$$
w=w_1\ldots w_l=u_1\ldots u_r,\;\;\;r,l\geq 1
$$
into the words $w_i, u_j$ ($1\leq i\leq l$, $1\leq j\leq r$)  from the sequence  $C$. In particular,  $w_1\neq u_1$. There exists the smallest number $i_0\geq 1$ such that $w_{i_0}=v_\mu$ or $u_{i_0}=v_\mu$ (otherwise  the words $w_i $ and $u_j$ all  belong to $C'$, which  implies that $C'$ is not a code). Without losing generality, we can assume that $w_{i_0}=v_\mu$. Then there exists $l'>l$ and the words $W_i $ ($1\leq i\leq l'$) in the sequence $C'$ such that $W_i=w_i$ for every $1\leq i<i_0$, $W_{i_0}=v_\kappa$, $W_{i_0+1}=v$ and the  equality $w_1\ldots w_l=W_1\ldots W_{l'}$ holds.

Suppose that $u_i\neq v_\mu$ for every $1\leq i\leq r$. Then every word $u_i$ belongs to $C'$. Since  $W_1\ldots W_{l'}=u_1\ldots u_r$ and  $C'$ is a code, we obtain $l'=r$ and $W_i=u_i$ for every $1\leq i\leq r$. Now, if $i_0>1$, then  $w_1=W_1=u_1$, contrary to our assumption. If $i_0=1$, then  $u_1=W_1=v_\kappa$, $u_2=W_2=v$ and hence the word $v=u_2$ belongs to $C$, contrary to our observation.

Thus, there exists  the smallest number $j_0\geq 1$ such that $u_{j_0}=v_\mu$. Then we have
$u_1\ldots u_r=U_1\ldots U_{r'}$ for some $r'>r$ and the words $U_i$ ($1\leq i\leq r'$) from the sequence $C'$ such that $U_i=u_i$ for every $1\leq i<j_0$, $U_{j_0}=v_\kappa$ and $U_{j_0+1}=v$. Since $w_1\ldots w_l=u_1\ldots u_r$, we obtain
$W_1\ldots W_{l'}=U_1\ldots U_{r'}$. Since $C'$ is a code, we obtain $l'=r'$ and $W_i=U_i$ for every $1\leq i\leq r'$. By the minimality of $i_0$, we have $j_0\geq i_0$. If $i_0>1$, then  $j_0\geq i_0>1$, and hence $u_1=U_1=W_1=w_1$, contrary to our assumption. Thus, it must be $i_0=1$, which implies $W_1=v_\kappa$, $W_2=v$. Now, if $j_0\geq 3$, then $u_2=U_2=W_2=v$, and hence the word $v=u_2$ would belong to $C$, which is impossible. If $j_0=2$, then $U_2=v_\kappa$ and, since $U_2=W_2=v$, we would obtain that the word $v=v_\kappa$ belongs to $C$. Thus it must be $j_0=i_0=1$. But then, by the definition of the numbers $i_0$, $j_0$, we have $w_1=u_1=v_\mu$. Thus the assumption that $C'$ is a code, always leads to a contradiction. Consequently, $C'$ is not a code, which finishes the proof of the lemma. \qed
\end{dwd}

Now, we are ready to prove our main result.
{
\renewcommand{\thetheorem}{\ref{t3}}
\begin{theorem}
For every $c\geq 2$ let us denote $L_c=(1,2,c)$. Then  we have:
$$
\lim_{c\to\infty}\rho_{2,L_c}=\frac{1}{6}=\alpha_2.
$$
\end{theorem}
\addtocounter{theorem}{-1}
}
\begin{dwd}
Let  $X=\{0,1\}$ be the binary alphabet. For every $c\geq 1$, we consider the subset $NUD(c)\subseteq X\times X^2\times X^c$ of those sequences which are not codes:
$$
NUD(c)=(X\times X^2\times X^c)\setminus UD_2(L_c).
$$
We have
\begin{equation}\label{udbin}
|UD_2(L_c)|=2^{c+3}-|NUD(c)|.
\end{equation}
For any letters $x,y,z\in X$, let us define the subset $K_{x,yz}(c)\subseteq NUD(c)$ as follows:
$$
K_{x,yz}(c):=\{(x,yz,w)\colon w\in X^c\}\subseteq NUD(c).
$$
Obviously, we have
$$
|K_{0,00}(c)|=|K_{1,11}(c)|=2^c,\;\;\;|K_{0,11}(c)|=|K_{1,00}(c)|.
$$
Since a sequence of words  is a code if and only if its reversal is a code (see~\cite{17,0} for example), we also have
$$
|K_{0,01}(c)|=|K_{0,10}(c)|=|K_{1,01}(c)|=|K_{1,10}(c)|.
$$
The set $NUD(c)$ is the union of the subsets $K_{x,yz}(c)$ ($x,y,z\in X$), and the following implication holds:
$$
(x,y,z)\neq (x',y',z')\Rightarrow K_{x,yz}(c)\cap K_{x',y'z'}(c)=\emptyset.
$$
Thus, we can write
\begin{eqnarray}\label{nud}
|NUD(c)|&=&2|K_{0,00}(c)|+2|K_{1,00}(c)|+4|K_{1,01}(c)|=\nonumber\\
&=&2^{c+1}+2|K_{1,00}(c)|+4|K_{1,01}(c)|.
\end{eqnarray}

Let $\widetilde{K}_{1,00}(c)$ be the set of all words $w\in X^c$ such that $(1,00,w)\in K_{1,00}(c)$, that is a word $w\in X^c$ belongs to $\widetilde{K}_{1,00}(c)$ if and only if the sequence $(1,00,w)$ is not a code. We also denote by $J_{1,00}(c)$ the set of all words $w\in X^c$  of the form
$$
1^{i_1}(00)^{j_1}\ldots 1^{i_k}(00)^{j_k},\;\;k\geq 1
$$
for some integers $i_l,j_l\geq 0$ ($1\leq l\leq k$).
\begin{stw}\label{lll3}
For every $c\geq 1$, we have $\widetilde{K}_{1,00}(c)=J_{1,00}(c)\cup\{0^c\}$. In particular, the following equalities hold
$$
|K_{1,00}(c)|=|\widetilde{K}_{1,00}(c)|=|J_{1,00}(c)\cup\{0^c\}|=F_{c+1}+(c)_2,
$$
where $F_n=\frac{(1+\sqrt{5})^n-(1-\sqrt{5})^n}{2^n\sqrt{5}}$  is the $n$-th Fibonacci number ($n=0,1,\ldots$) and  $(c)_2\in\{0,1\}$ is the remainder from the division of $c$ by $2$.
\end{stw}
\begin{dwd}[of Proposition~\ref{lll3}]
To show the equality $\widetilde{K}_{1,00}(c)=J_{1,00}(c)\cup\{0^c\}$, we use induction on $c$. We have $\widetilde{K}_{1,00}(1)=\{0,1\}=J_{1,00}(1)\cup\{0\}$. Suppose that there is $c\geq 1$ such that  $\widetilde{K}_{1,00}(i)=J_{1,00}(i)\cup\{0^i\}$ for every $1\leq i\leq c$. To derive the equality $\widetilde{K}_{1,00}(c+1)=J_{1,00}(c+1)\cup\{0^{c+1}\}$, we show at first that $J_{1,00}(c+1)\cup\{0^{c+1}\}$ is included in $\widetilde{K}_{1,00}(c+1)$. Indeed,  every word $w\in J_{1,00}(c+1)$ has at least two factorizations into the words $1$, $00$, $w$. For the word $w:=0^{c+1}$ we have $ww=(00)^{c+1}$, and hence $ww$ has two factorizations into the words $1$, $00$, $w$. Thus the set $J_{1,00}(c+1)\cup\{0^{c+1}\}$ is included in the set $\widetilde{K}_{1,00}(c+1)$.

To show the converse inclusion, let us choose $w\in \widetilde{K}_{1,00}(c+1)$ arbitrarily. Then one of the words $00$, $1$ must be a prefix of $w$ (otherwise the sequence $(1,00,w)$ would be a code). Thus there exists $v\in X^*$ such that $w=00v$ or $w=1v$, and hence, by Lemma~\ref{lll},  the sequence $(1,00,v)$ is not a code. Thus   $v\in \widetilde{K}_{1,00}(c-1)$ in the first case  and $v\in \widetilde{K}_{1,00}(c)$ in the second case. In the first case, by the inductive assumption, we obtain $v\in J_{1,00}(c-1)\cup\{0^{c-1}\}$ and  consequently
$$
w=00v\in J_{1,00}(c+1)\cup\{0^{c+1}\}.
$$
In the second case,  we have by the inductive assumption: $v\in J_{1,00}(c)\cup\{0^{c}\}$, and hence
$$
w=1v\in J_{1,00}(c+1)\cup\{10^{c}\}.
$$
Now, if $c$ is an even number, then $10^c\in J_{1,00}(c+1)$, and hence $w\in J_{1,00}(c+1)$. If $c$ is an odd number,  then  $w$ can not be equal to $10^c$. Indeed, if $c$ is odd, then the sequence $(1,00,10^c)$ is a code (to show this, we may use the Sardinas-Patterson algorithm, as then we have: $\mathcal{D}_0=\{1,00,10^c\}$, $\mathcal{D}_1=\{0^c\}$, and by the trivial induction on $i$, we have: $\mathcal{D}_i=\{0^{c-2i+2}\}$ for every $1\leq i\leq (c+1)/2$ and $\mathcal{D}_i=\{0\}$ for $i\geq (c+1)/2$. This implies $\mathcal{D}_i\cap \mathcal{D}_0=\emptyset$ for every $i\geq 1$). Thus, if $c$ is odd, then $w\neq 10^c$ and consequently $w\in J_{1,00}(c+1)$. Thus in every case, we have $w\in J_{1,00}(c+1)\cup\{0^{c+1}\}$ and hence   $\widetilde{K}_{1,00}(c+1)\subseteq J_{1,00}(c+1)\cup\{0^{c+1}\}$. The inductive argument finishes the proof of the first part.

To show the second part, its enough to show the equality
$$
|J_{1,00}(c)\cup\{0^c\}|=F_{c+1}+(c)_2.
$$
Let us denote $x_c:=|J_{1,00}(c)|$ for every $c\geq 1$.  Directly by the definition of the sets $J_{1,00}(c)$ ($c\geq 1$), we have for every $c\geq 3$:
\begin{eqnarray*}
1w\in J_{1,00}(c)&\Leftrightarrow& w\in J_{1,00}(c-1),\\
0w\in J_{1,00}(c)&\Leftrightarrow& w=0v\;\;{\rm and}\;\;v\in J_{1,00}(c-2).
\end{eqnarray*}
Thus $x_c=x_{c-1}+x_{c-2}$ for every $c\geq 3$. Since $x_1=1$ and $x_2=2$,  we obtain $x_c=F_{c+1}$ for every $c\geq 1$. Thus the equality $|J_{1,00}(c)|=F_{c+1}$ holds for every $c\geq 1$. Now, if $c$ is even, then $0^c\in J_{1,00}(c)$, and consequently
$$
|J_{1,00}(c)\cup\{0^c\}|=|J_{1,00}(c)|=F_{c+1}=F_{c+1}+(c)_2.
$$
If $c$ is odd, then $0^c\notin J_{1,00}(c)$, and consequently
$$
|J_{1,00}(c)\cup\{0^c\}|=|J_{1,00}(c)|+1=F_{c+1}+1=F_{c+1}+(c)_2.
$$
This finishes the proof of Proposition~\ref{lll3}
\qed
\end{dwd}

In the next step, we derive the number of elements of the set $K_{1,01}(c)$. We proceed in the similar way. Namely, we denote by $\widetilde{K}_{1,01}(c)$  the set of all words $w\in X^c$ such that $(1,01,w)\in K_{1,01}(c)$, that is a word $w\in X^c$ belongs to $\widetilde{K}_{1,01}(c)$ if and only if the sequence $(1,01,w)$ is not a code. We also denote by $J_{1,01}(c)$ the set of all words $w\in X^c$ which  do not contain two consecutive $0$'s.

\begin{stw}\label{lll6}
$\widetilde{K}_{1,01}(c)=J_{1,01}(c)$ for every $c\geq 1$. In particular, we have $|K_{1,01}(c)|=F_{c+2}$ for every $c\geq 1$.
\end{stw}
\begin{dwd}[of Proposition~\ref{lll6}]
The proof of the first part is by induction on $c$. For $c=1$, we have $\widetilde{K}_{1,01}(c)=\widetilde{K}_{1,01}(1)=\{0,1\}=J_{1,01}(1)=J_{1,01}(c)$. Suppose that there is $c\geq 1$ such that $\widetilde{K}_{1,01}(i)=J_{1,01}(i)$  for every $1\leq i\leq c$. To show the equality $\widetilde{K}_{1,01}(c+1)=J_{1,01}(c+1)$, we use (as before) the double inclusion argument. So let  $w\in  J_{1,01}(c+1)$ be arbitrary. Then there is $k\geq 1$ such that
$$
w=1^{i_1}01^{i_2}\ldots 01^{i_{k-1}}01^{i_k},
$$
where $i_1, i_k\geq 0$ and $i_t\geq 1$ for every $1<t<k$. In the case $i_k\geq 1$ we can write
$$
w=1^{i_1}(01)1^{i_2-1}\ldots (01)1^{i_k-1},
$$
which means that  $w$ has two different factorizations into the words $1$, $01$, $w$. If $i_k=0$, then we have
$$
w1=1^{i_1}(01)1^{i_2-1}\ldots (01)1^{i_{k-1}-1}01,
$$
which gives  that the word $w1$ has two different factorizations into the words $1$, $01$, $w$. Thus in every case, we have $w\in \widetilde{K}_{1,01}(c+1)$. Consequently, the inclusion $J_{1,01}(c+1)\subseteq\widetilde{K}_{1,01}(c+1)$ holds.

Conversely, let $w\in \widetilde{K}_{1,01}(c+1)$ be arbitrary. Then one of the words $01$, $1$ must be a prefix of $w$. Thus there exists $v\in X^*$ such that $w=01v$ or $w=1v$. By Lemma~\ref{lll}, we obtain that $(1,01,v)$ is not a code, which means that $v\in \widetilde{K}_{1,01}(c-1)$ in the first case and $v\in \widetilde{K}_{1,01}(c)$ in the second case. By the inductive assumption, the word $v$ does  not contain two consecutive $0$'s. But then the word $w$ also  does not contain two consecutive $0$'s, which means that $w\in J_{1,01}(c+1)$, and consequently  $\widetilde{K}_{1,01}(c+1)\subseteq J_{1,01}(c+1)$. The inductive argument finishes the proof of the first part.

To show the second part, it is enough to show the equality $|J_{1,01}(c)|=F_{c+2}$ for every $c\geq 1$. Let us  denote $y_c:=|J_{1,01}(c)|$. By the definition of the set $J_{1,01}(c)$, we  have for every $c\geq 3$:
\begin{eqnarray*}
1w\in J_{1,01}(c)&\Leftrightarrow& w\in J_{1,01}(c-1),\\
0w\in J_{1,01}(c)&\Leftrightarrow& w=1v\;\;{\rm and}\;\;v\in J_{1,01}(c-2).
\end{eqnarray*}
Thus $y_c=y_{c-1}+y_{c-2}$ for all $c\geq 3$. Since $y_1=2$ and $y_2=3$,  we obtain $y_c=F_{c+2}$ for every $c\geq 1$. This finishes the proof of Proposition~\ref{lll6}.\qed
\end{dwd}

Now, by the equalities~(\ref{udbin})--(\ref{nud}) and by Propositions~\ref{lll3}--\ref{lll6}, we can write:
\begin{eqnarray*}
|UD_2(L_c)|&=&2^{c+3}-(2^{c+1}+2|K_{1,00}(c)|+4|K_{1,01}(c)|)=\\
&=&2^{c+3}-2^{c+1}-2F_{c+1}-4F_{c+2}-2(c)_2=\\
&=&3\cdot2^{c+1}-2F_{c+4}-2(c)_2,\;\;\;c\geq 1,
\end{eqnarray*}
where the last equality follows from the identity
$$
2F_{c+1}+4F_{c+2}=2F_{c+4}.
$$
By Proposition~\ref{prefi}, we also have:
$$
|PR_2(L_c)|=|PR_2((1,2,c))|=2^c,\;\;\;c\geq 2.
$$
In consequence, we have for every $c\geq 2$:
$$
\rho_{2, L_c}=\frac{|PR_2(L_c)|}{|UD_2(L_c)|}=\frac{2^c}{3\cdot 2^{c+1}-2F_{c+4}-2(c)_2}=\frac{1}{6-\frac{F_{c+4}}{2^{c-1}}-\frac{(c)_2}{2^{c-1}}}.
$$
Since $\lim\limits_{c\to\infty}\frac{(c)_2}{2^{c-1}}=0$ and
\begin{eqnarray*}
\lim\limits_{c\to\infty}\frac{F_{c+4}}{2^{c-1}}=\lim\limits_{c\to\infty}\frac{(1+\sqrt{5})^{c+4}-(1-\sqrt{5})^{c+4}}{2^{2c+3}\sqrt{5}}=\\
=\frac{1}{2^{-5}\cdot\sqrt{5}}\lim\limits_{c\to\infty}\left(\left(\frac{1+\sqrt{5}}{4}\right)^{c+4}-\left(\frac{1-\sqrt{5}}{4}\right)^{c+4}\right)=0,
\end{eqnarray*}
we obtain $\lim\limits_{c\to\infty}\rho_{2, L_c}=\frac{1}{6}=\alpha_2$.\qed
\end{dwd}

\noindent
{\bf Acknowledgements.} We would like to thank to the anonymous referee for many valuable comments and accurate suggestions, which allow to simplify the presentation of the proof.


\begin{thebibliography}{00}

\bibitem{17} J. Berstel, D. Perrin, C. Reutenauer, {\it Codes and Automata (Encyclopedia of Mathematics and its Applications)} 1st Cambridge University Press New York, NY, USA, 2009 ISBN:052188831X 9780521888318

\bibitem {0} J. Berstel, D. Perrin, {\it Theory of codes} Pure and Applied Mathematics, vol. 117. Academic Press Inc., Orlando, FL, 1985.

\bibitem{1} E. K. Blum, {\it A note on free semigroups with two generators}, Bull. Amer. Math. Soc., vol. 71, pp. 678-679, 1965.

\bibitem{12} D. Derenourt, {\it  A three-word code which is not prefix-suffix composed}, Theoret. Comput. Sci., 163, 145-160, 1996.


\bibitem{15} J. Devolder, (1993). {\it Codes, mots infinis et bi-infinis}, These de doctorat, Universite Lille 1.

\bibitem{14} C.M. Fan, H.J. Shyr, S.S. Yu, {\it d-words and d-languages}, Acta Informatica 35:709–727, 1998.

\bibitem{2} J. Karhum\"{a}ki, {\it On three-element codes}, Theoret. Comput. Sci. 40: 3-11, 1985.

\bibitem{8} J. Karhum\"{a}ki, {\it  A property of three-element codes}, Theoret. Comput. Sci. 41: 215-222, 1985.

\bibitem{10} J. Karhum\"{a}ki, M. Latteux, I. Petre, {\it  The commutation with codes and ternary sets of words}. STACS, 2003.


\bibitem{16} M. Lothaire, {\it Algebraic combinatorics on words (Encyclopedia of Mathematics and its Applications 90)}, Cambridge University Press, 2002, ISBN 978-0-521-81220-7, MR 1905123.

\bibitem{9} Z. Z. Li, Y.S. Tsai, G.C.Yih, {\it  Characterizations on codes with three elements},  Soochow Journal of Mathematics 30(2): 177-196, 2004 (2004)

\bibitem{11} Z. Z. Li, Y.S. Tsai, {\it Three-element codes with one d-primitive word},  Acta Informatica, 41:171-180, 2004.

\bibitem{6} B. McMillan, {\it Two inequalities implied by unique decipherability},  IEEE Trans. Information Theory 2 (4): 115–116, (1956).

\bibitem{13} A. Restivo, S. Salemi, T. Sportelli, {\it Completing codes}, Theoret. Inform. Appl., 23, 135-147, 1989.

\bibitem{3} A. Sardinas, G. W. Patterson, {\it A necessary and sufficient condition for the unique decomposition of coded messages}, Convention Record of the I.R.E., 1953 National Convention, Part 8: Information Theory, pp. 104–108.

\bibitem{4} A. Woryna, {\it  On the set of uniquely decodable codes with a given sequence of code word lengths}, Discret. Math. 2017, vol. 340 iss. 2, s. 51-57.

\bibitem{22} A. Woryna, {\it On the ratio of prefix codes to all uniquely decodable codes with a given length distribution}, Discret. Appl. Math. 2018 vol. 244, s. 205-213.
\end{thebibliography}
\end{document}